# ASYMPTOTIC SPECTRAL THEORY FOR NONLINEAR TIME SERIES

By Xiaofeng Shao and Wei Biao Wu[1]

*University of Illinois at Urbana-Champaign and University of Chicago*

We consider asymptotic problems in spectral analysis of stationary causal processes. Limiting distributions of periodograms and smoothed periodogram spectral density estimates are obtained and applications to the spectral domain bootstrap are given. Instead of the commonly used strong mixing conditions, in our asymptotic spectral theory we impose conditions only involving (conditional) moments, which are easily verifiable for a variety of nonlinear time series.

**1. Introduction.** The frequency domain approach to time series analysis is an important subject; see [1, 7, 29] and [52] among others. An asymptotic distribution theory is needed, for example, in hypothesis testing and in the construction of confidence intervals. However, most of the asymptotic results developed in the literature are for strong mixing processes and processes with quite restrictive summability conditions on joint cumulants [6, 7, 56, 57]. Such conditions seem restrictive and they are not easily verifiable. For example, Andrews [2] showed that, for a simple autoregressive process with innovations being independent and identically distributed (i.i.d.) Bernoulli random variables, the process is not strong mixing. Other special processes discussed include Gaussian processes [60, 61] and linear processes [1].

There has been a recent surge of interest in nonlinear time series ([21, 53] and [65]). It seems that a systematic asymptotic spectral theory for such processes is lacking [11]. The primary goal of this paper is to establish an asymptotic spectral theory for stationary, causal processes. Let $(\varepsilon_n)_{n \in \mathbf{Z}}$ be a sequence of i.i.d. random variables; let

$$(1.1) \qquad X_n = G(\ldots, \varepsilon_{n-1}, \varepsilon_n),$$

Received January 2006.

[1]Supported in part by NSF Grant DMS-04-78704.

where $G$ is a measurable function such that $X_n$ is a proper random variable. Then the process $(X_n)$ is causal or nonanticipative in the sense that it only depends on $\mathcal{F}_n = (\ldots, \varepsilon_{n-1}, \varepsilon_n)$, not on the future innovations $\varepsilon_{n+1}, \varepsilon_{n+2}, \ldots$. The class of processes within the framework of (1.1) is quite large (cf. [53, 65, 66] and [74] among others).

Assume throughout this paper that $(X_n)_{n \in \mathbf{Z}}$ has mean zero and finite covariance function $r(k) = \mathbf{E}(X_0 X_k)$, $k \in \mathbf{Z}$. Let $i = \sqrt{-1}$ be the imaginary unit. If $(X_n)$ is short-range dependent, namely

$$(1.2) \qquad \sum_{k=0}^{\infty} |r(k)| < \infty,$$

then the spectral density

$$f(\lambda) = \frac{1}{2\pi} \sum_{k \in \mathbf{Z}} r(k) e^{ik\lambda}, \qquad \lambda \in \mathbf{R},$$

is continuous and bounded. Given the observations $X_1, \ldots, X_n$, let

$$S_n(\theta) = \sum_{k=1}^{n} X_k e^{ik\theta} \quad \text{and} \quad I_n(\theta) = \frac{1}{2\pi n} |S_n(\theta)|^2$$

be the Fourier transform and the periodogram, respectively. Let $\theta_k = 2\pi k/n$, $1 \le k \le n$, be the Fourier frequencies. Primary goals in spectral analysis include estimating the spectral density $f$ and deriving asymptotic distributions of $S_n(\theta)$ and $I_n(\theta)$.

We now introduce some notation. For a column vector $x = (x_1, \ldots, x_q)' \in \mathbf{R}^q$, let $|x| = (\sum_{j=1}^{q} x_j^2)^{1/2}$. Let $\xi$ be a random vector. Write $\xi \in \mathcal{L}^p$ ($p > 0$) if $\|\xi\|_p := [\mathbf{E}(|\xi|^p)]^{1/p} < \infty$ and let $\| \cdot \| = \| \cdot \|_2$. For $\xi \in \mathcal{L}^1$ define projection operators $\mathcal{P}_k \xi = \mathbf{E}(\xi | \mathcal{F}_k) - \mathbf{E}(\xi | \mathcal{F}_{k-1})$, $k \in \mathbf{Z}$, where we recall $\mathcal{F}_k = (\ldots, \varepsilon_{k-1}, \varepsilon_k)$. For two positive sequences $(a_n)$, $(b_n)$, denote by $a_n \asymp b_n$ that there exists a constant $c$ such that $0 < c \le a_n/b_n \le 1/c < \infty$ for all large $n$ and by $a_n \sim b_n$ that $a_n/b_n \to 1$ as $n \to \infty$. Let $C > 0$ denote a generic constant which may vary from line to line; let $\Phi$ be the standard normal distribution function. Denote by "$\Rightarrow$" convergence in distribution and by $N(\mu, \sigma^2)$ a normal distribution with mean $\mu$ and variance $\sigma^2$. All asymptotic statements in the paper are with respect to $n \to \infty$ unless otherwise specified.

The paper is structured as follows. In Section 2 we shall establish a central limit theorem for the Fourier transform $S_n(\theta)$ at Fourier frequencies. Asymptotic properties of smoothed periodogram estimates of $f$ are discussed in Section 3. Section 4 shows the consistency of the frequency domain bootstrap approximation to sampling distributions of spectral density estimates for both linear and nonlinear processes. Section 5 gives sufficient conditions for geometric moment contraction [see (3.1)], a basic dependence assumption used in this paper. Some examples are also presented in that section. Proofs are gathered in the Appendix.

**2. Fourier transforms.** The periodogram is a fundamental quantity in frequency domain analysis. Its asymptotic analysis has a substantial history; see, for example, [57], Theorem 5.3, page 131, for mixing processes; [8], Theorem 10.3.2, page 347, [63] and [70] for linear processes. Other contributions can be found in [38, 46, 55, 71] and [77]. Recently, in a general setting, Wu [73] considered asymptotic distributions of $S_n(\theta)$ at a fixed $\theta$. However, results in [73] do not apply to $S_n(\theta)$ at the Fourier frequencies. Here we shall show that $S_n(\theta_k)$ are asymptotically independent normals under mild conditions; see Theorem 2.1 below. The central limit theorem is applied to empirical distribution functions of normalized periodogram ordinates (cf. Corollary 2.2). In the literature the latter problem has been mainly studied for i.i.d. random variables [25, 26, 36] and linear processes [12].

Denote the real and imaginary parts of $S_n(\theta_j)/\sqrt{\pi n f(\theta_j)}$ by

$$Z_j = \frac{\sum_{k=1}^n X_k \cos(k\theta_j)}{\sqrt{\pi n f(\theta_j)}}, \qquad Z_{j+m} = \frac{\sum_{k=1}^n X_k \sin(k\theta_j)}{\sqrt{\pi n f(\theta_j)}}, \qquad j = 1, \ldots, m,$$

where $m = m_n := \lfloor (n-1)/2 \rfloor$ and $\lfloor a \rfloor$ is the integer part of $a$. Let $\Omega_p = \{c \in \mathbf{R}^p : |c| = 1\}$ be the unit sphere. For the set $J = \{j_1, \ldots, j_p\}$ with $1 \le j_1 < \cdots < j_p \le 2m$ write the vector $Z_J = (Z_{j_1}, \ldots, Z_{j_p})'$. Let the class $\Xi_{m,p} = \{J \subset \{1, \ldots, 2m\} : \#J = p\}$, where $\#J$ is the cardinality of $J$.

THEOREM 2.1. *Assume* $X_t \in \mathcal{L}^2$,

$$(2.1) \qquad \kappa := \sum_{k=0}^\infty \|\mathcal{P}_0 X_k\| < \infty$$

*and* $f_* := \min_{\theta \in \mathbf{R}} f(\theta) > 0$. *Then for any fixed* $p \in \mathbf{N}$, *we have*

$$\sup_{J \in \Xi_{m,p}} \sup_{c \in \Omega_p} \sup_{x \in \mathbf{R}} |\mathbf{P}(Z_J' c \le x) - \Phi(x)| = o(1) \qquad \text{as } n \to \infty.$$

Theorem 2.1 asserts that the projection of any vector of $p$ of the $Z_j$'s on any direction is asymptotically normal. The condition (2.1) was first proposed by Hannan [30]. In many situations it is easily verifiable since it only involves conditional moments. For generalizations see [75]. In the special case of linear processes $X_t = \sum_{j=0}^\infty a_j \varepsilon_{t-j}$, where $\varepsilon_j$ are i.i.d. with mean 0 and finite variance and $\sum_{j=0}^\infty a_j^2 < \infty$, (2.1) becomes $\sum_{j=0}^\infty |a_j| < \infty$, indicating that $(X_n)$ is short-range dependent. In the literature, central limit theorems are established for Fourier transforms of linear processes ([21], page 63; [8], page 347, among others). The spectral density may be unbounded if (2.1) is violated.

COROLLARY 2.1. *Let $q \in \mathbf{N}$. Under the conditions of Theorem* 2.1, *we have*

$$\left\{ \frac{S_n(\theta_{l_j})}{\sqrt{n\pi f(\theta_{l_j})}}, 1 \le j \le q \right\} \Rightarrow \{Y_{2j-1} + iY_{2j}, 1 \le j \le q\}$$

*for integers $1 \le l_1 < l_2 < \cdots < l_q \le m$, where the indices $l_j$ may depend on $n$, and $Y_k$, $1 \le k \le 2q$, are i.i.d. standard normals. Consequently, for $\tilde{I}_n(\theta) := I_n(\theta)/f(\theta)$,*

$$\{\tilde{I}_n(\theta_{l_j}), 1 \le j \le q\} \Rightarrow \{E_j, 1 \le j \le q\},$$

*where $E_j$ are i.i.d. standard exponential random variables* [exp(1)].

Corollary 2.1 easily follows from Theorem 2.1 via the Cramér–Wold device. Let

$$F_{\tilde{I},m}(x) := \frac{1}{m} \sum_{j=1}^{m} \mathbf{1}_{\tilde{I}_n(\theta_j) \le x}$$

be the empirical distribution function of $\tilde{I}_n(\theta_k)$ and $F_E(x) := 1 - e^{-x}$, $x \ge 0$.

COROLLARY 2.2. *Under the conditions of Theorem* 2.1, *we have*

$$(2.2) \qquad \sup_{x \ge 0} |F_{\tilde{I},m}(x) - F_E(x)| \to 0 \quad \text{in probability.}$$

PROOF. Since $F_{\tilde{I},m}$ and $F_E$ are nondecreasing, it suffices to show (2.2) for a *fixed* $x$. Let $p_j = p_j(x) = \mathbf{P}[\tilde{I}_n(\theta_j) \le x]$ and $p_{j,k} = p_{j,k}(x) = \mathbf{P}[\tilde{I}_n(\theta_j) \le x, \tilde{I}_n(\theta_k) \le x]$; let $U$ and $V$, independent of the process $(X_j)$, be i.i.d. uniformly distributed over $\{1, \ldots, m\}$. By Corollary 2.1, $p_U \to F_E(x)$ and $p_{U,V} \to F_E(x)^2$ almost surely. By the Lebesgue dominated convergence theorem, $\mathbf{E}(p_U) \to F_E(x)$ and $\mathbf{E}(p_{U,V}) \to F_E(x)^2$. Notice that

$$\mathbf{E}(p_U) = m^{-1} \sum_{j=1}^{m} p_j \quad \text{and} \quad \mathbf{E}(p_{U,V}) = m^{-2} \sum_{j=1}^{m} \sum_{k=1}^{m} p_{j,k}.$$

So $\|F_{\tilde{I},m}(x) - F_E(x)\|^2 = \mathbf{E}(p_{U,V}) - F_E^2(x) + 2F_E(x)\{F_E(x) - \mathbf{E}(p_U)\}$ and (2.2) follows. □

REMARK 2.1. The above argument also implies that, for any integer $k \ge 2$,

$$\sup_{x_1, \ldots, x_k \ge 0} \left| m^{-k} \sum_{j_1=1}^{m} \cdots \sum_{j_k=1}^{m} \mathbf{1}_{\tilde{I}_n(\theta_{j_1}) \le x_1, \ldots, \tilde{I}_n(\theta_{j_k}) \le x_k} - \prod_{j=1}^{k} F_E(x_j) \right| \to 0$$

in probability.

Fay and Soulier [22] obtained a functional central limit theorem for $F_{\tilde{I},m}(x)$ for i.i.d. random variables. It seems very difficult to generalize their results to the nonlinear case.

**3. Spectral density estimation.** Given a realization $(X_j)_{j=1}^n$, the spectral density $f$ can be estimated by

$$f_n(\lambda) = \int_{-\pi}^{\pi} W_n(\lambda - \mu) I_n(\mu) \, d\mu,$$

where $W_n(\lambda)$ is a smoothing weight function [cf. (3.2)]. Here we study asymptotic properties of the smoothed periodogram estimate $f_n$. Spectral density estimation is an important problem and there is a rich literature. However, restrictive structural conditions have been imposed in many earlier results. For example, Brillinger [6] assumed that all moments exist and cumulants of all orders are summable. Anderson [1] dealt with linear processes. Rosenblatt [56] considered strong mixing processes and assumed the summability condition of cumulants up to the eighth order. Due to those limitations, the classical results cannot be directly applied to nonlinear time series. Recently, Chanda [11] obtained asymptotic normality of $f_n$ for a class of nonlinear processes. However, it seems that his formulation does not include popular nonlinear time series models including GARCH, EXPAR and ARMA–GARCH; see Section 5 for examples.

To establish an asymptotic theory for $f_n$, we shall adopt the geometric-moment contraction (GMC) condition. Let $(\varepsilon_k')_{k \in \mathbf{Z}}$ be an i.i.d. copy of $(\varepsilon_k)_{k \in \mathbf{Z}}$; let $X_n' = G(\ldots, \varepsilon_{-1}', \varepsilon_0', \varepsilon_1, \ldots, \varepsilon_n)$ be a coupled version of $X_n$. We say that $X_n$ is GMC($\alpha$), $\alpha > 0$, if there exist $C > 0$ and $0 < \rho = \rho(\alpha) < 1$ such that, for all $n \in \mathbf{N}$,

$$(3.1) \qquad \mathbf{E}(|X_n' - X_n|^\alpha) \le C\rho^n.$$

Inequality (3.1) indicates that the process $(X_n)$ quickly "forgets" the past $\mathcal{F}_0 = (\ldots, \varepsilon_{-1}, \varepsilon_0)$. Note that under GMC(2), $|r(k)| = O(\rho^k)$ for some $\rho \in (0, 1)$ and hence the spectral density function is infinitely many times differentiable.

Many nonlinear time series models satisfy GMC (cf. Section 5). Moreover, the GMC condition provides a convenient framework for a limit theory for nonlinear time series; see [32, 75] and [76]. In view of those features, instead of the widely used strong mixing condition, we employ the GMC as an underlying assumption for our asymptotic theory of spectral density estimates.

Let $\hat{r}(k) = n^{-1} \sum_{j=1}^{n-|k|} X_j X_{j+|k|}$, $|k| < n$, be the estimated covariances; let $a(\cdot)$ be an even, Lipschitz continuous function with support $[-1, 1]$ and $a(0) = 1$; let $B_n$ be a sequence of positive integers with $B_n \to \infty$ and $B_n/n \to$

$0$; let $b_n = 1/B_n$,

$$W_n(\lambda) = \frac{1}{2\pi} \sum_{k=-B_n}^{B_n} a(kb_n) e^{-ik\lambda} \quad \text{and}$$

(3.2)

$$f_n(\lambda) = \frac{1}{2\pi} \sum_{k=-B_n}^{B_n} \hat{r}(k) a(kb_n) e^{-ik\lambda}.$$

THEOREM 3.1. *Assume* (3.1), $X_n \in \mathcal{L}^{4+\delta}$ *for some* $\delta > 0$, $B_n \to \infty$ *and* $B_n = o[n/(\log n)^{2+8/\delta}]$. *Then*

(3.3)                $\sqrt{nb_n}\{f_n(\lambda) - \mathbf{E}(f_n(\lambda))\} \Rightarrow N(0, \sigma^2(\lambda)),$

*where* $\sigma^2 := \sigma^2(\lambda) = \{1 + \eta(2\lambda)\} f^2(\lambda) \int_{-1}^{1} a^2(t)\, dt$ *and* $\eta(\lambda) = 1$ *if* $\lambda = 2k\pi$ *for some integer* $k$ *and* $\eta(\lambda) = 0$ *otherwise.*

REMARK 3.1. The GMC has this interesting property: If $X_n \in \mathcal{L}^p, p > 0$, and GMC($\alpha_0$) holds for some $\alpha_0 > 0$, then $X_n$ is GMC($\alpha$) for all $\alpha \in (0, p)$ ([75], Lemma 2).

By Remark 3.1, the moment condition $X_n \in \mathcal{L}^{4+\delta}$ in Theorem 3.1 together with GMC($\alpha$) implies GMC(4) and consequently the absolute summability of cumulants up to the fourth order (cf. Lemmas A.1 and A.2). In the context of strong mixing processes, Rosenblatt ([57], page 138) imposed $X_n \in \mathcal{L}^8$. Rosenblatt [57] also posed the problem of whether the eighth-order cumulant summability condition can be weakened to fourth order. Theorem 3.1 partially solves the conjecture for nonlinear processes satisfying GMC under the moment condition $X_n \in \mathcal{L}^{4+\delta}$. Additionally, Theorem 3.1 is applicable to a variety of nonlinear time series models (Section 5) that are not covered by Chanda [11].

Joint asymptotic distributions of spectral density estimates at different frequencies (cf. Corollary 3.1 below) follow from the arguments in [48], Theorem 5A and [56] since GMC(4) ensures the summability of the fourth cumulants; see Lemma A.2.

COROLLARY 3.1. *Let* $\lambda_1, \ldots, \lambda_s \in [0, \pi]$ *be* $s$ *different frequencies. Then under the conditions of Theorem* 3.1, $\sqrt{nb_n}\{f_n(\lambda_j) - \mathbf{E}(f_n(\lambda_j))\}$, $j = 1, \ldots, s$, *are jointly asymptotically independent* $N(0, \sigma^2(\lambda_j))$, $j = 1, \ldots, s$.

The problem of maximum deviation of spectral density estimates has been studied by Woodroofe and Van Ness [72] for linear processes and Rudzkis [58] for Gaussian processes. For nonlinear processes, we have:

THEOREM 3.2. *Assume* (3.1), $X_n \in \mathcal{L}^{4+\delta}$ *for some* $\delta \in (0, 4]$, $B_n \to \infty$, $B_n = O(n^\eta)$, $0 < \eta < \delta/(4 + \delta)$ *and* $f_* := \min_{\mathbf{R}} f(\theta) > 0$. *Then*

$$(3.4) \qquad \max_{\lambda \in [0,\pi]} \sqrt{nb_n} |f_n(\lambda) - \mathbf{E}(f_n(\lambda))| = O_{\mathbf{P}}((\log n)^{1/2}).$$

Under GMC(2), since $\|\mathcal{P}_0 X_k\| = O(\rho^k)$, we have (2.1). However, it is quite difficult to establish (3.3) under the weaker condition (2.1). Regarding (3.4), for linear processes the distributional result in [72] implies that the bound $O_{\mathbf{P}}((\log n)^{1/2})$ is optimal. We are unable to obtain a similar distributional result for nonlinear processes.

For long memory processes, (1.2) is violated and $f$ may not be well defined, so Theorems 3.1 and 3.2 are not applicable. A simple example is the fractionally integrated process $(1 - B)^d X_j = \varepsilon_j$, where $0 < d < 1/2$ is the long memory parameter, $B$ is the back-shift operator and $\varepsilon_j$ are i.i.d. with mean 0 and finite variance. Then the spectral density $f(\lambda) \asymp |\lambda|^{-2d}$ as $\lambda \to 0$ and $f(0)$ is not well defined. In this case an important problem is to estimate $d$; see [54] and [59] and references cited therein.

## 4. Frequency domain bootstrap.
Here we consider bootstrap approximations of the distribution of the lag window estimate (3.2). Bootstrapping in the frequency domain has received considerable attention. See [33, 45] and [64] for Gaussian processes and [24, 37] and [47] for linear processes. For nonlinear processes we adopt the residual-based bootstrap procedure proposed by Franke and Härdle [24]. A variant of it is discussed in Remark 4.4. Let $I_j = I(\omega_j)$, $\omega_j = 2\pi j/n$, $j \in F_n = \{-\lfloor (n-1)/2 \rfloor, \ldots, \lfloor n/2 \rfloor\}$. Note that $\hat{r}(k) = n^{-1} 2\pi \sum_{j \in F_n} I_j e^{ik\omega_j}$. Then the lag window estimate (3.2) can be written as

$$f_n(\lambda) = \frac{1}{2\pi} \sum_{k=-B_n}^{B_n} \hat{r}(k) a(kb_n) e^{-ik\lambda} = \frac{1}{n} \sum_{j \in F_n} I_j \sum_{k=-B_n}^{B_n} a(kb_n) e^{-ik(\lambda - \omega_j)}.$$

(4.1)
The bootstrap procedure consists of the following several steps:

1. Calculate periodogram ordinates $\{I_j\}$, $j = 1, \ldots, N := \lfloor n/2 \rfloor$.

2. Obtain an estimate $\tilde{f}$ of $f$ (e.g., a lag window estimate with bandwidth $\tilde{b}_n := \tilde{B}_n^{-1}$).

3. Let $\bar{\varepsilon}_j = \tilde{\varepsilon}_j / \bar{\varepsilon}$, where $\tilde{\varepsilon}_j = I_j / \tilde{f}_j$, $\tilde{f}_j = \tilde{f}(\omega_j)$ and $\bar{\varepsilon} = N^{-1} \sum_{j=1}^N \tilde{\varepsilon}_j$.

4. Draw i.i.d. bootstrap samples $\{\varepsilon_j^*\}$ from the empirical distribution of $\bar{\varepsilon}_j$.

5. Let $I_j^* = \tilde{f}_j \varepsilon_j^*$ be the bootstrapped periodograms; let $I_{-j}^* = I_j^*$ and $I_0^* = 0$.

The rescaling treatment in step 3 avoids an unpleasant bias at the resampling stage. Setting $I_0^* = 0$ in step 5 corresponds to the fact that, for a mean-corrected sample, the periodogram value is 0 at frequency 0. The sampling distribution of $g_n(\lambda) = \sqrt{nb_n}\{f_n(\lambda) - f(\lambda)\}$ is expected to be close to its bootstrap counterpart $g_n^*(\lambda) = \sqrt{nb_n}\{f_n^*(\lambda) - \tilde{f}(\lambda)\}$, where

$$f_n^*(\lambda) = \frac{1}{n} \sum_{j \in F_n} I_j^* \sum_{k=-B_n}^{B_n} a(kb_n) e^{-ik(\lambda - \omega_j)}$$

is the bootstrapped version of (4.1). Here we measure the closeness by Mallows' $d_2$ metric [4]. For two probability measures $P_1$ and $P_2$ on $\mathbf{R}$ with $\int_{\mathbf{R}} |x|^2 \, dP_j < \infty$, $j = 1, 2$, let $d_2(P_1, P_2) = \inf \|Y_1 - Y_2\|$, where the infimum is taken over all vectors $(Y_1, Y_2)$ with marginal distributions $P_1$ and $P_2$. Write

$$d_2[g_n(\lambda), g_n^*(\lambda)] = d_2\{\mathbf{P}[g_n(\lambda) \in \cdot], \mathbf{P}[g_n^*(\lambda) \in \cdot \,|X_1, \ldots, X_n]\}.$$

The bootstrap procedure is said to be (weakly) consistent if $d_2[g_n(\lambda), g_n^*(\lambda)] = o_{\mathbf{P}}(1)$. Let $\mathcal{L}(\cdot \,|X_1, \ldots, X_n)$ denote the conditional distribution given the sample $X_1, \ldots, X_n$.

It seems that in the literature the theoretical investigation of the consistency problem has been limited to linear processes. Let $X_t = \sum_{j=-\infty}^{\infty} a_j \varepsilon_{t-j}$. Franke and Härdle [24] proved the consistency of their residual-based procedure under the condition

(4.2)        $\sup\{|\mathbf{E}(e^{iu\varepsilon_1})|; |u| \geq \delta\} < 1$        for all $\delta > 0$.

Condition (4.2) excludes many interesting cases. For example, it is violated if $\varepsilon_1$ is a Bernoulli random variable. Franke and Härdle [24] conjectured that their results still hold without (4.2). The latter condition is removed in Corollary 4.1 of Theorem 4.1 below at the expense of the stronger eighth moment condition. Theorem 4.1 is also applicable to nonlinear processes; see Corollary 4.2. Since our results hold under various combinations of conditions, it is convenient to label the common ones:

(A1)  $\lim_{x \to 0} x^{-2}\{1 - a(x)\} = c_2$, where $c_2$ is a nonzero constant.
(A2)  $\min_{\lambda \in [0, \pi]} f(\lambda) > 0$.
(A3)  $\max_{\lambda \in [0, \pi]} |\tilde{f}(\lambda) - f(\lambda)| = o_{\mathbf{P}}(b_n)$.
(A3′)  $\max_{\lambda \in [0, \pi]} |\tilde{f}(\lambda) - f(\lambda)| = o_{\mathbf{P}}(1)$.
(A4)  $\sum_{k \in \mathbf{Z}} |r(k)| k^2 < \infty$.
(A4′)  $\sum_{k \in \mathbf{Z}} |r(k)| k < \infty$.
(A5)  $\sum_{t_1, \ldots, t_{k-1} \in \mathbf{Z}} |\mathrm{cum}(X_0, X_{t_1}, \ldots, X_{t_{k-1}})| < \infty$ for $k = 3, 4$.
(A5′)  $\sum_{t_1, \ldots, t_{k-1} \in \mathbf{Z}} |\mathrm{cum}(X_0, X_{t_1}, \ldots, X_{t_{k-1}})| < \infty$ for $k = 3, \ldots, 8$.
(A6)  $\sqrt{nb_n}\{f_n(\lambda) - \mathbf{E}(f_n(\lambda))\} \Rightarrow N(0, \sigma^2(\lambda))$ and $nb_n \, \mathrm{var}(f_n(\lambda)) \to \sigma^2(\lambda)$.

REMARK 4.1. Condition (A1) says that $a(\cdot)$ is locally quadratic at 0 and it is satisfied for many lag windows. It is related to the bias. By Anderson [1], Theorem 9.4.3, or Priestley [52], page 459, under (A1), (A4) and $B_n^3 = o(n)$,

$$(4.3) \qquad B_n^2\{\mathbf{E}(f_n(\lambda)) - f(\lambda)\} \to c_2 f''(\lambda),$$

$$\text{where } f''(\lambda) = -\frac{1}{2\pi}\sum_{k\in\mathbf{Z}} r(k)k^2 e^{-ik\lambda}.$$

Additionally, if (A6) holds, then the optimal bandwidth $b_n$ is of order $n^{-1/5}$ in the sense of mean square error.

REMARK 4.2. The cumulant summability conditions (A5) and (A5') are commonly imposed in spectral analysis [7, 57]. For the linear process $X_t = \sum_{j=-\infty}^{\infty} a_j \varepsilon_{t-j}$ with $\sum_{j=-\infty}^{\infty}|a_j| < \infty$, (A5) [resp. (A5')] holds if $\varepsilon_1 \in \mathcal{L}^4$ [resp. $\varepsilon_1 \in \mathcal{L}^8$]. By Lemma A.1, for the process (1.1), (A5) [resp. (A5')] is satisfied under GMC(4) [resp. GMC(8)]. Zhurbenko and Zuev [79] and Andrews [3] considered strong mixing processes.

Let $\mathbf{P}^*$, $\mathbf{E}^*$ and var$^*$ denote the conditional probability, expectation and variance given $X_j, 1 \leq j \leq n$; let $V_n(\lambda) = \sqrt{nb_n}\{f_n(\lambda) - \mathbf{E}(f_n(\lambda))\}$, $V_n^*(\lambda) = \sqrt{nb_n}\{f_n^*(\lambda) - \mathbf{E}^* f_n^*(\lambda)\}$, $\beta_n(\lambda) = \sqrt{nb_n}\{\mathbf{E}(f_n(\lambda)) - f(\lambda)\}$ and $\beta_n^*(\lambda) = \sqrt{nb_n} \times \{\mathbf{E}^* f_n^*(\lambda) - \tilde{f}(\lambda)\}$. For the consistency of the bootstrap approximation, it is common to treat the variance and the bias separately.

PROPOSITION 4.1. Assume $X_t \in \mathcal{L}^8$, (A2), (A3), (A4'), (A5') and (A6). Let $B_n^2 = o(n)$. Then $d_2[V_n(\lambda), V_n^*(\lambda)] \to 0$ in probability.

PROPOSITION 4.2. Assume $X_t \in \mathcal{L}^4$, (A1), (A4) and (A5). Let $b_n = o(\tilde{b}_n)$, $B_n^3 = o(n)$ and $\tilde{B}_n^5 = o(n)$. Then $B_n^2\{\mathbf{E}^* f_n^*(\lambda) - \tilde{f}(\lambda)\} \to c_2 f''(\lambda)$ in probability.

REMARK 4.3. The condition $b_n = o(\tilde{b}_n)$ is needed to ensure the consistency of the bias; see (4.3). So $\tilde{f}(\lambda)$ is smoother than $f_n(\lambda)$. Oversmoothing is a common practice in the frequency domain bootstrap [24, 37, 47].

THEOREM 4.1. Assume $X_t \in \mathcal{L}^8$, (A1), (A4), (A5') and (A6). Let $b_n \asymp n^{-1/5}$ and $b_n = o(\tilde{b}_n)$. Then $d_2[g_n(\lambda), g_n^*(\lambda)] = o_{\mathbf{P}}(1)$ and $d_2[g_n(\lambda)/f(\lambda), g_n^*(\lambda)/\tilde{f}(\lambda)] = o_{\mathbf{P}}(1)$.

PROOF. In the proof $\lambda$ is suppressed and we write $g_n$, and so on, for $g_n(\lambda)$, and so on. Since $d_2^2(g_n, g_n^*) = d_2^2(V_n, V_n^*) + d_2^2(\beta_n, \beta_n^*)$ ([4], Lemma 8.8), by Propositions 4.1, 4.2 and (4.3), $d_2(g_n, g_n^*) = o_{\mathbf{P}}(1)$. The second assertion follows similarly. By (A2), (A3) and Proposition 4.2, $\beta_n^*/\tilde{f} - \beta_n/f = (\beta_n^* - \beta_n)/\tilde{f} + (\tilde{f}^{-1} - f^{-1})\beta_n = o_{\mathbf{P}}(1)$. It remains to verify $d_2(V_n/f, V_n^*/\tilde{f}) = o_{\mathbf{P}}(1)$. By Lemma 8.3 in [4], it suffices in view of (A6) to show that $\mathrm{var}^*(V_n^*/\tilde{f}) \to \sigma^2/f^2$ and $\mathcal{L}(V_n^*/\tilde{f}|X_1,\ldots,X_n) \Rightarrow N(0,\sigma^2/f^2)$ in probability. By (A2) and (A3), these two assertions follow from relation (A.21) in the proof of Proposition 4.1. □

REMARK 4.4. Since the residuals $\{I_n(\omega_j)/f(\omega_j)\}$ are asymptotically i.i.d. exp(1) (Corollary 2.1), a modified procedure is to replace the bootstrapped residuals $\varepsilon_j^*$ by i.i.d. standard exponential variables. For this modified bootstrap procedure, Theorem 4.1 holds with the assumption (A5′) replaced by (A5) and the eighth moment condition weakened to $X_t \in \mathcal{L}^4$; see the proof of Proposition 4.1.

COROLLARY 4.1. *Let* $X_t = \sum_{j=-\infty}^{\infty} a_j \varepsilon_{t-j}$, *where* $|a_k| = O(|k|^{-1-\beta})$, $\beta > 1/5$ *and* $\varepsilon_1 \in \mathcal{L}^8$. *Assume* (A1), (A2), (A4), $b_n \asymp n^{-1/5}$ *and* $\tilde{b}_n \asymp n^{-\eta_1}$, $\eta_1 \in (1/10, 1/5)$. *Then the conclusions in Theorem* 4.1 *hold.*

PROOF. By Theorem 4.1, it suffices to verify (A3), (A5′) and (A6). (A6) follows from Theorems 9.3.4 and 9.4.1 in [1]. The assumption (A5′) is satisfied under $\mathbf{E}(\varepsilon_1^8) < \infty$ and $|a_k| = O(|k|^{-1-\beta})$, $\beta > 1/5$ (see Remark 4.2). Note that

$$
\begin{aligned}
(4.4) \quad & \max_{\lambda \in [0,\pi]} |\tilde{f}(\lambda) - f(\lambda)| \\
& \leq \max_{\lambda \in [0,\pi]} |\tilde{f}(\lambda) - \mathbf{E}(\tilde{f}(\lambda))| + \max_{\lambda \in [0,\pi]} |\mathbf{E}(\tilde{f}(\lambda)) - f(\lambda)|,
\end{aligned}
$$

which is of order $O_{\mathbf{P}}((\log n)^{1/2}/(n\tilde{b}_n)^{1/2}) + O_{\mathbf{P}}(\tilde{b}_n^2) = o_{\mathbf{P}}(b_n)$ by Theorem 2.1 in [72] and (4.3). So (A3) follows. □

COROLLARY 4.2. *Let the process* (1.1) *satisfy GMC*(8). *Assume* (A1), (A2), $b_n \asymp n^{-1/5}$ *and* $\tilde{b}_n \asymp n^{-\eta_2}$, $\eta_2 \in (1/10, 1/5)$. *Then the conclusions in Theorem* 4.1 *hold.*

PROOF. We shall apply Theorem 4.1. By Lemma A.1, GMC(8) implies (A4) and (A5′), while (A6) [resp. (A3)] follows from Theorem 3.1 [resp. Theorem 3.2 and (4.4)]. □

5. **Applications.** There are two popular criteria to check the stationarity of nonlinear time series models, drift-type conditions [10, 23, 42, 67, 68, 69] and contraction conditions [16, 19, 34, 76]. It turns out that contraction conditions typically imply GMC under some extra mild assumptions, and are thus quite useful in proving limit theorems [32, 75]. In this section we consider nonlinear autoregressive models and present sufficient conditions for GMC so that our asymptotic spectral theory is applicable.

Let $\varepsilon, \varepsilon_n$ be i.i.d., $p, d \geq 1$; let $X_n \in \mathbf{R}^d$ be recursively defined by

$$(5.1) \qquad X_{n+1} = R(X_n, \ldots, X_{n-p+1}; \varepsilon_{n+1}),$$

where $R$ is a measurable function. Suitable conditions on $R$ implies GMC.

THEOREM 5.1. *Let $\alpha > 0$ and $\alpha' = \min(1, \alpha)$. Assume that $R(y_0; \varepsilon) \in \mathcal{L}^\alpha$ for some $y_0$ and that there exist constants $a_1, \ldots, a_p \geq 0$ such that $\sum_{j=1}^p a_j < 1$ and*

$$(5.2) \qquad \|R(y; \varepsilon) - R(y'; \varepsilon)\|_\alpha^{\alpha'} \leq \sum_{j=1}^p a_j |x_j - x_j'|^{\alpha'}$$

*holds for all $y = (x_1, \ldots, x_p)$ and $y' = (x_1', \ldots, x_p')$. Then* (i) (5.1) *admits a stationary solution of the form* (1.1) *and* (ii) $X_n$ *satisfies GMC($\alpha$). In particular, if there exist functions $H_j$ such that $|R(y; \varepsilon) - R(y'; \varepsilon)| \leq \sum_{j=1}^p H_j(\varepsilon)|x_j - x_j'|$ for all $y$ and $y'$ and $\sum_{j=1}^p \|H_j(\varepsilon)\|_\alpha^{\alpha'} < 1$, then we can let $a_j = \|H_j(\varepsilon)\|_\alpha^{\alpha'}$.*

Duflo [18] assumed $\alpha \geq 1$ and called (5.2) the Lipschitz mixing condition. We allow $\alpha < 1$. Similar conditions are given in [27].

PROOF OF THEOREM 5.1. It follows from the arguments in [76] and Lemma 6.2.10 and Proposition 6.3.22 in [18]. For completeness we include the proof here. Without loss of generality let $d = 1$. Let $\alpha < 1$. For $y = (x_1, \ldots, x_p) \in \mathbf{R}^p$ define the random map $R_\varepsilon(y) = (R(y, \varepsilon), x_1, \ldots, x_{p-1})$. Let $Z_m(y)$ be the first element of the vector $R_{\varepsilon_0} \circ R_{\varepsilon_{-1}} \circ \cdots \circ R_{\varepsilon_{-m}}(y)$, where $m$ is a nonnegative integer. By (5.2), we have for $m \geq p$ that

$$\|Z_m(y) - Z_m(y')\|_\alpha^\alpha \leq \sum_{j=1}^p a_j \|Z_{m-j}(y) - Z_{m-j}(y')\|_\alpha^\alpha.$$

Since $a_1, \ldots, a_p$ are nonnegative and $\sum_{j=1}^p a_j < 1$, it is easily seen that the preceding relation implies that there exist constants $C > 0$ and $\lambda_0 \in (0, 1)$ depending only on $a_1, \ldots, a_p$ and $\alpha$ such that

$$(5.3) \qquad \|Z_m(y) - Z_m(y')\|_\alpha^\alpha \leq C\lambda_0^m |y - y'|^\alpha$$

holds for all $m \geq 0$. See also Lemma 6.2.10 in [18]. Applying (5.3) with $y = y_0$ and $y' = R_{\varepsilon_{-m-1}}(y_0)$, since $\lambda_0 < 1$ and $\alpha < 1$,

$$\mathbf{E}\left(\sum_{m=0}^\infty |Z_m(y_0) - Z_{m+1}(y_0)|\right)^\alpha \leq \sum_{m=0}^\infty \|Z_m(y_0) - Z_{m+1}(y_0)\|_\alpha^\alpha$$

$$\leq C\sum_{m=0}^\infty \lambda_0^m \|y_0 - R_\varepsilon(y_0)\|_\alpha^\alpha < \infty.$$

So $\{Z_m(y_0)\}_{m \geq 0}$ is a Cauchy sequence and it has an almost sure limit $Z_\infty$ (say) which is in $\mathcal{L}^\alpha$. Since $Z_\infty$ is $\mathcal{F}_0$-measurable, we can write $Z_\infty = G(\mathcal{F}_0)$ for some measurable function $G$. By (5.3), for any $y$, $Z_m(y)$ converges almost surely to the same limit $Z_\infty$. So we can express $X_n = G(\mathcal{F}_n)$, $n \in \mathbf{Z}$. Let $\mathcal{F}_j^* = (\ldots, \varepsilon_{j-1}', \varepsilon_j')$. By stationarity, (ii) follows from (5.3) by letting $y = (G(\mathcal{F}_{-m-1}), \ldots, G(\mathcal{F}_{-m-p}))$ and $y' = (G(\mathcal{F}_{-m-1}^*), \ldots, G(\mathcal{F}_{-m-p}^*))$. The other case $\alpha \geq 1$ can be similarly dealt with. See Proposition 6.3.22 in [18]. $\square$

THEOREM 5.2. *Let $(\eta_t)$ satisfy GMC$(\alpha)$; let $\theta_1, \ldots, \theta_p, \phi_1, \ldots, \phi_q$, $p, q \in \mathbf{N}$, be real coefficients and the roots of the equation $\lambda^p - \sum_{k=1}^p \theta_k \lambda^{p-k} = 0$ lie inside the unit circle. Then the autoregressive moving average (ARMA) $(p, q)$ process $X_t$ defined below also satisfies GMC$(\alpha)$:*

$$X_t - \theta_1 X_{t-1} - \cdots - \theta_p X_{t-p} = \eta_t - \phi_1 \eta_{t-1} - \cdots - \phi_q \eta_{t-q}.$$

Theorem 5.2 shows that the GMC property is preserved in ARMA modeling [43] and that it is an easy consequence of the representation $X_t = \sum_{k=0}^\infty b_k \eta_{t-k}$ with $|b_k| \leq C\rho^k$ for some $\rho \in (0, 1)$. Min [43] considered the case $\alpha \geq 1$. Theorem 5.2 implies that the ARMA–ARCH and ARMA–GARCH models [39] are GMC; see Examples 5.4 and 5.5.

Near-epoch dependence (NED) is widely used in econometrics for central limit theorems [14, 15]. The process (1.1) is geometrically NED [G-NED$(\alpha)$] on $(\varepsilon_s)$ in $L_\alpha$, $\alpha > 0$, if there exist $C < \infty$ and $\rho \in (0, 1)$ such that

$$\|X_t - \mathbf{E}(X_t | \varepsilon_{t-m}, \varepsilon_{t-m+1}, \ldots, \varepsilon_t)\|_\alpha \leq C\rho^m$$

holds for all $m \in \mathbf{N}$. It is easily seen that, for $\alpha \geq 1$, GMC$(\alpha)$ is equivalent to G-NED$(\alpha)$. In certain situations GMC is more convenient to work with; see Remark 5.1. Additionally, GMC has the nice property that $X_t'$ is identically distributed as $X_t$, while in NED the distribution of $\mathbf{E}(X_t | \varepsilon_{t-m}, \ldots, \varepsilon_t)$ typically differs. Here we list some examples that are not covered by Davidson [15].

EXAMPLE 5.1. Amplitude-dependent exponential autoregressive (EXPAR) models have been studied by Jones [35]. Let $\varepsilon_j \in \mathcal{L}^\alpha$ be i.i.d. innovations and

$$X_n = [\alpha_1 + \beta_1 \exp(-aX_{n-1}^2)]X_{n-1} + \varepsilon_n,$$

where $\alpha_1, \beta_1, a > 0$ are real parameters. Then $H_1(\varepsilon) = |\alpha_1| + |\beta_1|$. By Theorem 5.1, $X_n$ is GMC$(\alpha)$ if $|\alpha_1| + |\beta_1| < 1$.

EXAMPLE 5.2. Let $\theta_1, \ldots, \theta_5$ be real parameters and consider the AR(2) model with ARCH(2) errors [20],

$$X_n = \theta_1 X_{n-1} + \theta_2 X_{n-2} + \varepsilon_n \sqrt{\theta_3^2 + \theta_4^2 X_{n-1}^2 + \theta_5^2 X_{n-2}^2}.$$

Theorem 5.1 is applicable: we can let $H_1(\varepsilon) = |\theta_1| + |\varepsilon \theta_4|$ and $H_2(\varepsilon) = |\theta_2| + |\varepsilon \theta_5|$. Then GMC($\alpha$), $\alpha > 0$, holds if $\sum_{j=1}^{2} \|H_j(\varepsilon)\|_\alpha^\alpha < 1$ and $\varepsilon_1 \in \mathcal{L}^\alpha$.

Let $A_t$ be $p \times p$ random matrices and $B_t$ be $p \times 1$ random vectors. The generalized random coefficient autoregressive process $(X_t)$ is defined by

$$(5.4) \qquad X_{t+1} = A_{t+1} X_t + B_{t+1}, \qquad t \in \mathbf{Z}.$$

Let $(A_t, B_t)$ be i.i.d. Bilinear and GARCH models fall within the framework of (5.4). The stationarity, geometric ergodicity and $\beta$-mixing properties have been studied by Pham [50], Mokkadem [44] and Carrasco and Chen [9]. Their results require that innovations have a density, which is not needed in our setting.

For a $p \times p$ matrix $A$, let $|A|_\alpha = \sup_{z \neq 0} |Az|_\alpha / |z|_\alpha$, $\alpha \geq 1$, be the matrix norm induced by the vector norm $|z|_\alpha = (\sum_{j=1}^{p} |z_j|^\alpha)^{1/\alpha}$. It is easily seen that $X_t$ is GMC($\alpha$), $\alpha \geq 1$, if $\mathbf{E}(|A_0|_\alpha) < 1$ and $\mathbf{E}(|B_0|_\alpha) < \infty$. By Jensen's inequality, we have $\mathbf{E}(\log|A_0|_\alpha) < 0$. By Theorem 1.1 of [5],

$$(5.5) \qquad X_n = \sum_{k=0}^{\infty} A_n A_{n-1} \cdots A_{n-k+1} B_{n-k}$$

converges almost surely.

EXAMPLE 5.3. Consider the subdiagonal bilinear model [28, 62]

$$(5.6) \qquad X_t = \sum_{j=1}^{p} a_j X_{t-j} + \sum_{j=0}^{q} c_j \varepsilon_{t-j} + \sum_{j=0}^{P} \sum_{k=1}^{Q} b_{jk} X_{t-j-k} \varepsilon_{t-k}.$$

Let $s = \max(p, P + q, P + Q)$, $r = s - \max(q, Q)$ and $a_{p+j} = 0 = c_{q+j} = b_{P+k,Q+j} = 0$, $k, j \geq 1$; let $H$ be a $1 \times s$ vector with the $(r+1)$st element 1 and all others 0, $c$ be an $s \times 1$ vector with the first $r-1$ elements 0 followed by 1, $a_1 + c_1, \ldots, a_{s-r} + c_{s-r}$, and $d$ be an $s \times 1$ vector with the first $r$ elements 0 followed by $b_{01}, \ldots, b_{0,s-r}$. Define the $s \times s$ matrices

$$A = \begin{pmatrix} 0 & 1 & & 0 & & 0 \\ & & \ddots & & 0 & \\ 0 & & & 1 & & 0 \\ 0 & 0 & & a_1 & \ddots & 0 \\ & & & & \vdots & 1 \\ a_s & \cdots & \cdots & a_{s-r} & & 0 \end{pmatrix},$$

$$B = \begin{pmatrix} 0 & \cdots & 0 & 0 & \cdots & 0 \\ \vdots & \vdots & \vdots & \vdots & \vdots & \vdots \\ 0 & \cdots & 0 & 0 & \cdots & 0 \\ b_{r1} & \cdots & b_{01} & 0 & \cdots & 0 \\ \vdots & \vdots & \vdots & \vdots & \vdots & \vdots \\ b_{r,s-r} & \cdots & b_{0,s-r} & 0 & \cdots & 0 \end{pmatrix}.$$

Let $Z_t$ be an $s \times 1$ vector with the $j$th entry $X_{t-r+j}$ if $1 \le j \le r$ and

$$\sum_{k=j}^{r} a_k X_{t+j-k} + \sum_{k=j}^{s-r} \left\{ c_k + \sum_{l=0}^{P} b_{lk} X_{t+j-k-l} \right\} \varepsilon_{t+j-k}$$

if $1 + r \le j \le s$. Pham [49, 51] gave the representation

$$(5.7) \qquad X_t = H Z_{t-1} + \varepsilon_t, \qquad Z_t = (A + B\varepsilon_t) Z_{t-1} + c\varepsilon_t + d\varepsilon_t^2.$$

By (5.7), $X_t$ is GMC($\alpha$), $\alpha \ge 1$, if $\varepsilon_1 \in \mathcal{L}^{2\alpha}$ and $\mathbf{E}(|A + B\varepsilon_1|^\alpha) < 1$. By (5.5), $Z_t$ admits a casual representation and so does $X_t$.

REMARK 5.1. Davidson [15] considered the bilinear model (5.6) with $q = 0$ and $Q = 1$. He commented that, due to the complexity of moment expressions, it is not easy to show G-NED(2) for general cases. In comparison, our argument works.

EXAMPLE 5.4. Ding, Granger and Engle [17] proposed the asymmetric GARCH($r, s$) model

$$(5.8) \quad X_t = \varepsilon_t \sqrt{h_t}, \qquad h_t^{\varsigma/2} = \alpha_0 + \sum_{j=1}^{r} \alpha_j (|X_{t-j}| - \gamma X_{t-j})^\varsigma + \sum_{j=1}^{s} \beta_j h_{t-j}^{\varsigma/2},$$

where $\alpha_0 > 0$, $\alpha_j \ge 0$ $(j = 1, \ldots, r)$ with at least one $\alpha_j > 0$, $\beta_j \ge 0$ $(j = 1, \ldots, s)$, $\varsigma \ge 0$ and $|\gamma| < 1$. The linear GARCH($r, s$) model is a special case of (5.8) with $\varsigma = 2, \gamma = 0$. Wu and Min [75] showed GMC for linear GARCH($r, s$) models. Let $Z_t = (|\varepsilon_t| - \gamma\varepsilon_t)^\varsigma, \xi_{\varsigma t} = (\alpha_0 Z_t, 0, \ldots, \alpha_0, 0, \ldots, 0)'_{(r+s) \times 1}$, of which the $(r+1)$st element is $\alpha_0$ and

$$A_{\varsigma t} = \left( \begin{array}{ccc|ccc} \alpha_1 Z_t & \cdots & \alpha_r Z_t & \beta_1 Z_t & \cdots & \beta_s Z_t \\ & I_{(r-1) \times (r-1)} & O_{(r-1) \times 1} & & O_{(r-1) \times s} & \\ \hline \alpha_1 & \cdots & \alpha_r & \beta_1 & \cdots & \beta_s \\ & O_{(s-1) \times r} & & & I_{(s-1) \times (s-1)} & O_{(s-1) \times 1} \end{array} \right).$$

Ling and McAleer [40] showed that $X_t \in \mathcal{L}^{m\varsigma}$ for some $m \in \mathbf{N}$ if and only if

$$(5.9) \qquad \Delta\{\mathbf{E}(A_{\varsigma t}^{\otimes m})\} < 1,$$

where $\otimes$ is the usual Kronecker product and $\Delta(A)$ is the largest eigenvalue of the matrix $(A'A)^{1/2}$. Further, $X_t$ admits a casual representation (1.1); see

Theorem 3.1 of Ling and McAleer [40]. It turns out that (5.9) also implies $GMC(m\varsigma)$.

PROPOSITION 5.1. *For the asymmetric GARCH$(r,s)$ model* (5.8), *let* $\varepsilon_t \in \mathcal{L}^{m\varsigma}$, $\varsigma \geq 1$. *Then* $X_t$ *is* $GMC(m\varsigma)$ *if* (5.9) *holds.*

PROOF. Let $Y_t = [(|X_t| - \gamma X_t)^\varsigma, \ldots, (|X_{t-r+1}| - \gamma X_{t-r+1})^\varsigma, h_t^{\varsigma/2}, \ldots, h_{t-s}^{\varsigma/2}]'$. Then $Y_t = A_{\varsigma t} Y_{t-1} + \xi_{\varsigma t}$ [40]. Let $Y_0'$, independent of $\{\varepsilon_t, t \in \mathbf{Z}\}$, be an i.i.d. copy of $Y_0$. We recursively define $Y_t' = A_{\varsigma t} Y_{t-1}' + \xi_{\varsigma t}$, $t \geq 1$. Let $\tilde{Y}_t = Y_t - Y_t'$. Then $\tilde{Y}_t = A_{\varsigma t} \tilde{Y}_{t-1}$. Applying the argument of Proposition 3 in [75], we have

$$\tilde{Y}_t^{\otimes m} = A_{\varsigma t}^{\otimes m} \tilde{Y}_{t-1}^{\otimes m} = \cdots = A_{\varsigma t}^{\otimes m} \cdots A_{\varsigma 1}^{\otimes m} \tilde{Y}_0^{\otimes m}.$$

Thus $\mathbf{E}(\tilde{Y}_t^{\otimes m}) = [\mathbf{E}(A_{\varsigma 1}^{\otimes m})]^t \mathbf{E}(\tilde{Y}_0^{\otimes m})$ since $A_{\varsigma t}, \ldots, A_{\varsigma 1}$ are i.i.d. By (5.9), $|\mathbf{E}(\tilde{Y}_t^{\otimes m})| \leq C\rho^t$ for some $\rho \in (0,1)$. So $\mathbf{E}(|h_t^{\varsigma/2} - (h_t')^{\varsigma/2}|^m)$ is bounded by $C\rho^t$ and

$$\mathbf{E}(|X_t - X_t'|^{m\varsigma}) = \mathbf{E}(\varepsilon_t^{m\varsigma}) \mathbf{E}(|\sqrt{h_t} - \sqrt{h_t'}|^{m\varsigma}) \leq C\mathbf{E}(|h_t^{\varsigma/2} - (h_t')^{\varsigma/2}|^m) \leq C\rho^t,$$

where the inequality $|a - b|^\varsigma \leq |a^\varsigma - b^\varsigma|$, $a \geq 0, b \geq 0, \varsigma \geq 1$, is applied. □

EXAMPLE 5.5. Let $\varepsilon_t$ be i.i.d. with mean 0 and variance 1. Consider the signed volatility model [78]

$$(5.10) \qquad X_t = \varepsilon_t |s_t|^{1/\varsigma}, \qquad s_t = g(\varepsilon_{t-1}) + c(\varepsilon_{t-1}) s_{t-1}, \qquad \varsigma > 0.$$

When $s_t = h_t^\varsigma > 0$, (5.10) reduces to the general GARCH$(1,1)$ model [31, 41]

$$X_t = \varepsilon_t h_t, \qquad h_t^\varsigma = g(\varepsilon_{t-1}) + c(\varepsilon_{t-1}) h_{t-1}^\varsigma, \qquad \varsigma > 0.$$

We shall show that the model (5.10) satisfies GMC under mild conditions.

PROPOSITION 5.2. *For the model* (5.10), *let* $\mathbf{E}(|\varepsilon_1|^{\alpha\varsigma}) < 1$, $\mathbf{E}\{|c(\varepsilon_1)|^\alpha\} < 1$ *and* $g(\varepsilon_1) \in \mathcal{L}^\alpha$, $\alpha > 0$. *Let* $\varsigma \geq 1$. *Then* $X_t$ *is* $GMC(\varsigma\alpha)$.

PROOF. By Theorem 5.1, $s_t$ is $GMC(\alpha)$. Since $\mathbf{E}\{(|s_t|^{1/\varsigma} - |s_t'|^{1/\varsigma})^{\varsigma\alpha}\} \leq \mathbf{E}(|s_t - s_t'|^\alpha)$ and $X_t = \varepsilon_t |s_t|^{1/\varsigma}$, $X_t$ is $GMC(\varsigma\alpha)$. □

Since $\mathbf{E}\{|c(\varepsilon_1)|^\alpha\} < 1$ implies $\mathbf{E}\{\log|c(\varepsilon_1)|\} < 0$, by Theorem 1 of [78], $X_t$ has a unique stationary solution which admits the casual representation (1.1).

## APPENDIX

We now give the proofs of the results in Sections 2–4.

**A.1. Proof of Theorem 2.1.** For presentational clarity we restrict $J = \{j_1, \ldots, j_p\} \subset \{1, \ldots, m\}$ and hence $Z_{j_l}$ corresponds to the real parts of $S_n(\theta_{j_l})$. The argument easily extends to general cases. Let

$$T_n = \sum_{k=1}^n \mu_k X_k, \qquad \text{where } \mu_k = \mu_k(c, J) = \sum_{l=1}^p \frac{c_l \cos(k\theta_{j_l})}{\sqrt{\pi f(\theta_{j_l})}}, \qquad 1 \le k \le n.$$

Since $f_* := \min_{\mathbf{R}} f(\theta) > 0$, there exists $\mu_*$ such that $|\mu_k| \le \mu_*$ for all $c \in \Omega_p$ and $J \in \Xi_{m,p}$. Let $d_n(h) = n^{-1} \sum_{k=1+h}^n \mu_k \mu_{k-h}$ if $0 \le h \le n-1$ and $d_n(h) = 0$ if $h \ge n$. Note that

$$\sum_{k=1}^n \cos(k\theta_{j_l}) \cos[(k+h)\theta_{j_{l'}}] = \frac{n}{2} \cos(h\theta_{j_l}) \mathbf{1}_{j_l = j_{l'}}.$$

Then it is easily seen that there exists a constant $K_0 > 0$ such that for all $h \ge 0$,

$$\tau_n(h) = \sup_{J \in \Xi_{m,p}} \sup_{c \in \Omega_p} \left| d_n(h) - \sum_{l=1}^p c_l^2 \frac{\cos(h\theta_{j_l})}{2\pi f(\theta_{j_l})} \right| \le \frac{K_0 h}{n}.$$

Clearly $\tau_n(h) \le \mu_* + (2\pi f_*)^{-1} =: K_1$. So we have uniformly over $J$ and $c$ that

$$
\text{(A.1)}
\begin{aligned}
\left| \frac{\|T_n\|^2}{n} - 1 \right| &= \left| d_n(0) r(0) + 2 \sum_{h=1}^\infty d_n(h) r(h) - 1 \right| \\
&\le 2 \sum_{h=0}^\infty \tau_n(h) r(h) \le \sum_{h=0}^\infty K_2 \min(h/n, 1) r(h) \to_{n \to \infty} 0
\end{aligned}
$$

by the Lebesgue dominated convergence theorem, where $K_2 = 2(K_0 + K_1)$.

Let $\tilde{T}_n = \sum_{k=1}^n \mu_k \tilde{X}_k$, where $\tilde{X}_k = \mathbf{E}(X_k | \varepsilon_{k-\ell+1}, \ldots, \varepsilon_k)$ are $\ell$-dependent. So $\delta_\ell = \|X_0 - \tilde{X}_0\| \to 0$ as $\ell \to \infty$. If $k < \ell$, then $\mathcal{P}_0 \tilde{X}_k = \mathbf{E}(\mathcal{P}_0 X_k | \varepsilon_{k-\ell+1}, \ldots, \varepsilon_0)$. By Jensen's inequality $\|\mathcal{P}_0 \tilde{X}_k\| \le \|\mathcal{P}_0 X_k\|$. If $k \ge \ell$, then $\mathcal{P}_0 \tilde{X}_k = 0$. Clearly $\|\mathcal{P}_0(X_k - \tilde{X}_k)\| \le 2\delta_\ell$. By the Lebesgue dominated convergence theorem, (2.1) entails that

$$
\text{(A.2)}
\begin{aligned}
\frac{\|T_n - \tilde{T}_n\|}{\sqrt{n}} &= \left[ \frac{1}{n} \sum_{j=-\infty}^n \|\mathcal{P}_j(T_n - \tilde{T}_n)\|^2 \right]^{1/2} \le \mu_* \sum_{k=0}^\infty \|\mathcal{P}_0(X_k - \tilde{X}_k)\| \\
&\le \mu_* \sum_{k=0}^\infty 2 \min(\|\mathcal{P}_0 X_k\|, \delta_\ell) \to_{\ell \to \infty} 0.
\end{aligned}
$$

Let $g_n(r) = r^2 \mathbf{E}[\tilde{X}^2 \mathbf{1}(|\tilde{X}| \ge \sqrt{n}/r)]$. Since $\mathbf{E}(\tilde{X}^2) < \infty$, $\lim_{n \to \infty} g_n(r) = 0$ for any fixed $r > 0$. Note that $g_n$ is nondecreasing in $r$. Then there exists a sequence $r_n \uparrow \infty$ such that $g_n(r_n) \to 0$. Let $Y_k = \tilde{X}_k \mathbf{1}(|\tilde{X}_k| \le \sqrt{n}/r_n)$ and $T_{n,Y} = \sum_{k=1}^{n} \mu_k Y_k$. Then $\|Y_k - \tilde{X}_k\| = o(1/r_n)$. Since $Y_k - \tilde{X}_k$ are $\ell$-dependent,

$$\text{(A.3)} \qquad \|T_{n,Y} - \tilde{T}_n\| \le \sum_{a=1}^{\ell} \left\| \sum_{b \le n, \ell | (b-a)} \mu_b (Y_b - \tilde{X}_b) \right\| = o(\sqrt{n}/r_n),$$

where $\ell | h$ means that $\ell$ is a divisor of $h$. Let $p_n = \lfloor r_n^{1/4} \rfloor$ and blocks $B_t = \{a \in \mathbf{N} : 1 + (t-1)(p_n + \ell) \le a \le p_n + (t-1)(p_n + \ell)\}$, $1 \le t \le t_n := \lfloor 1 + (n - p_n)/(p_n + \ell) \rfloor$. Define $U_t = \sum_{a \in B_t} \mu_a Y_a$, $V_n = \sum_{t=1}^{t_n} U_t$, $R_n = T_{n,Y} - V_n$, $W = (V_n - \mathbf{E}(V_n))/\sqrt{n}$ and $\Delta = \tilde{T}_n/\sqrt{n} - W$. Then $U_t$ are independent and $\|R_n\| = O(\sqrt{t_n})$ since $Y_a$ are $\ell$-dependent. Note that $|\mathbf{E}(V_n)| = O(n)|\mathbf{E}(Y_k)| = o(\sqrt{n}/r_n)$. Then by (A.3),

$$\text{(A.4)} \qquad \begin{aligned} \sqrt{n}\|\Delta\| &\le |\mathbf{E}(V_n)| + \|V_n - \tilde{T}_n\| = o(\sqrt{n}/r_n) + O(\sqrt{t_n} + \sqrt{n}/r_n) \\ &= O(\sqrt{t_n}). \end{aligned}$$

Since $|U_t|^3 \le \mu_*^3 p_n^2 \sum_{a \in B_t} |Y_a|^3$ and $\mathbf{E}(Y_a^2) \le \mathbf{E}(X_k^2)$, $\mathbf{E}(|U_t|^3) = O(p_n^3 \sqrt{n}/r_n)$. By the Berry–Esseen theorem ([13], page 304),

$$\text{(A.5)} \qquad \begin{aligned} \sup_x |\mathbf{P}(W \le x) - \Phi(x/\|W\|)| &\le C \sum_{t=1}^{t_n} \mathbf{E}(|U_t|^3) \times \|V_n - \mathbf{E}(V_n)\|^{-3} \\ &= O(t_n p_n^3 \sqrt{n}/r_n) \times n^{-3/2} = O(p_n^{-2}). \end{aligned}$$

Let $\delta = \delta_n = p_n^{-1/4}$. By (A.4), (A.5) and

$$\text{(A.6)} \qquad \begin{aligned} \mathbf{P}(W \le w - \delta) - \mathbf{P}(|\Delta| \ge \delta) &\le \mathbf{P}(W + \Delta \le w) \\ &\le \mathbf{P}(W \le w + \delta) + \mathbf{P}(|\Delta| \ge \delta), \end{aligned}$$

we have $\sup_x |\mathbf{P}(\tilde{T}_n \le \sqrt{n}x) - \Phi(\sqrt{n}x/\|\tilde{T}_n\|)| = O[p_n^{-2} + \mathbf{P}(|\Delta| \ge \delta) + \delta + \delta^2] = O(\delta)$ since $\sup_x |\Phi(x/\sigma_1) - \Phi(x/\sigma_2)| \le C|\sigma_1/\sigma_2 - 1|$ holds for some constant $C$.

Let $W_1 = \tilde{T}_n/\sqrt{n}$, $\Delta_1 = (T_n - \tilde{T}_n)/\sqrt{n}$ and $\eta = \eta_{\ell,n} = (\|T_n - \tilde{T}_n\|/\sqrt{n})^{1/2}$. We apply (A.6) with $W, \Delta$ replaced by $W_1, \Delta_1$,

$$\sup_x \left| \mathbf{P}\left( \frac{T_n}{\sqrt{n}} \le x \right) - \Phi\left( \frac{\sqrt{n}x}{\|T_n\|} \right) \right| = O(\mathbf{P}(|\Delta_1| \ge \eta) + \delta + \eta + \eta^2).$$

So the conclusion follows from (A.1) and (A.2) by first letting $n \to \infty$ and then $\ell \to \infty$.

**A.2. Proof of Theorem 3.1.** The following two lemmas are needed.

LEMMA A.1 ([76]). *Assume* (3.1) *with* $\alpha = k$ *for some* $k \in \mathbf{N}$. *Then there exists a constant* $C > 0$ *such that for all* $0 \le m_1 \le \cdots \le m_{k-1}$,

$$|\operatorname{cum}(X_0, X_{m_1}, \ldots, X_{m_{k-1}})| \le C\rho^{m_{k-1}/[k(k-1)]}.$$

LEMMA A.2. *Let the sequence* $s_n \in \mathbf{N}$ *satisfy* $s_n \le n$ *and* $B_n = o(s_n)$; *let*

$$(A.7) \qquad Y_u := Y_u(\lambda) = (2\pi)^{-1} \sum_{k=-B_n}^{B_n} X_u X_{u+k} a(kb_n) \cos(k\lambda).$$

*Then under GMC*(4) *we have* $\|\sum_{u=1}^{s_n} \{Y_u - \mathbf{E}(Y_u)\}\|^2 \sim s_n B_n \sigma^2$.

PROOF. Let $L(s) = \{(m_1, m_2, m_3) \in \mathbf{Z}^3 : \max_{1 \le i \le 3} |m_i| = s\}$ and $c(m_1, m_2, m_3) = \operatorname{cum}(X_0, X_{m_1}, X_{m_2}, X_{m_3})$. So $\#L(s) \le 6(2s+1)^2$. By Lemma A.1, we have

$$\sum_{m_1, m_2, m_3 \in \mathbf{Z}} |c(m_1, m_2, m_3)| \le C \sum_{s=0}^{\infty} \sum_{(m_1, m_2, m_3) \in L(s)} |c(m_1, m_2, m_3)|$$

$$\le C \sum_{s=0}^{\infty} s^2 \rho^{s/[4(4-1)]} < \infty.$$

See also Remark 3 in [76]. Then the lemma follows from equations (3.9)–(3.12) in [56], page 1174. □

PROOF OF THEOREM 3.1. Let $\rho = \rho(4)$, $\alpha_k = a(kb_n)\cos(k\lambda)$ and

$$h_n(\lambda) := \frac{1}{2\pi\sqrt{nB_n}} \left( \sum_{k=0}^{B_n} \sum_{u=n-k+1}^{n} X_u X_{u+k} \alpha_k + \sum_{k=-B_n}^{-1} \sum_{u=n+k+1}^{n} X_u X_{u+k} \alpha_k \right).$$

By the summability of cumulants of orders 2 and 4 (cf. [57], page 139), $\|h_n(\lambda)\| = (nB_n)^{-1/2}O(B_n)$. Recall (A.7) for the definition of $Y_u$ and let $g_n := g_n(\lambda) = \sum_{u=1}^{n} Y_u(\lambda)$. Then

$$(A.8) \qquad \sqrt{nb_n}\{f_n(\lambda) - \mathbf{E}(f_n(\lambda))\} = \frac{g_n - \mathbf{E}(g_n)}{\sqrt{nB_n}} + h_n(\lambda) - \mathbf{E}(h_n(\lambda)).$$

For $k \in \mathbf{Z}$ let $\tilde{X}_k = \mathbf{E}(X_k | \varepsilon_{k-l+1}, \ldots, \varepsilon_k)$, where $l = l_n = \lfloor c \log n \rfloor$ and $c = -8/\log \rho$. Let $\tilde{Y}_u := \tilde{Y}_u(\lambda)$ be the corresponding sum with $X_k$ replaced by $\tilde{X}_k$. Observe that $\tilde{X}_n$ and $\tilde{X}_m$ are i.i.d. if $|n - m| \ge l$ and $\tilde{Y}_u$ and $\tilde{Y}_v$ are i.i.d. if $|u - v| \ge 2B_n + l$. The independence plays an important role in establishing the asymptotic normality of $\tilde{g}_n := \tilde{g}_n(\lambda) = \sum_{u=1}^n \tilde{Y}_u(\lambda)$. Then $\|g_n - \tilde{g}_n\| = o(1)$ since

$$(A.9) \quad \|Y_u - \tilde{Y}_u\| \le (2\pi)^{-1} \sum_{k=-B_n}^{B_n} \|X_u X_{u+k} - \tilde{X}_u \tilde{X}_{u+k}\| |\alpha_k| = O(B_n \rho^{l/4}).$$

Let $\psi_n = n/(\log n)^{2+8/\delta}$, $p_n = \lfloor \psi_n^{2/3} B_n^{1/3} \rfloor$ and $q_n = \lfloor \psi_n^{1/3} B_n^{2/3} \rfloor$. Then

$$(A.10) \quad \begin{aligned} &p_n, q_n \to \infty, \qquad q_n = o(p_n), \\ &2B_n + l = o(q_n) \quad \text{and} \quad k_n = \lfloor n/(p_n + q_n) \rfloor \to \infty. \end{aligned}$$

Define the blocks $\mathcal{L}_r = \{j \in \mathbf{N} : (r-1)(p_n + q_n) + 1 \le j \le r(q_n + p_n) - q_n\}$, $1 \le r \le k_n$, $\mathcal{S}_r = \{j \in \mathbf{N} : r(p_n + q_n) - q_n + 1 \le j \le r(q_n + p_n)\}$, $1 \le r \le k_n - 1$ and $\mathcal{S}_{k_n} = \{j \in \mathbf{N} : k_n(p_n + q_n) - q_n + 1 \le j \le n\}$. Let $U_r = \sum_{j \in \mathcal{L}_r} \tilde{Y}_j$ and $V_r = \sum_{j \in \mathcal{S}_r} \tilde{Y}_j$. Observe that $U_1, \ldots, U_{k_n}$ are i.i.d. and $V_1, \ldots, V_{k_n-1}$ are also i.i.d. By Lemma A.2 and (A.9),

$$(A.11) \quad \begin{aligned} \|U_1 - \mathbf{E}(U_1)\| &= \left\| \sum_{j=1}^{p_n} \{Y_j - \mathbf{E}(Y_0)\} \right\| + O(p_n \|Y_0 - \tilde{Y}_0\|) \\ &\sim (p_n B_n \sigma^2)^{1/2} + O(p_n B_n \rho^{l/4}) \sim (p_n B_n \sigma^2)^{1/2}. \end{aligned}$$

Similarly, $\|V_1 - \mathbf{E}(V_1)\| \sim (q_n B_n \sigma^2)^{1/2} + O(q_n B_n \rho^{l/4})$. By (A.10),

$$\begin{aligned} \mathrm{var}(V_1 + \cdots + V_{k_n}) &= (k_n - 1)\|V_1 - \mathbf{E}(V_1)\|^2 + \|V_{k_n} - \mathbf{E}(V_{k_n})\|^2 \\ &= O(k_n q_n B_n) + O[(p_n + q_n) B_n] = o(n B_n). \end{aligned}$$

Then we have $(n B_n)^{-1/2} \{\tilde{g}_n - \mathbf{E}(\tilde{g}_n)\} \Rightarrow N(0, \sigma^2)$ if

$$(A.12) \quad (n B_n)^{-1/2} \sum_{r=1}^{k_n} \{U_r - \mathbf{E}(U_1)\} \Rightarrow N(0, \sigma^2).$$

Let $\tau = 2 + \delta/2$. Case (i) $[\log n = o(B_n)]$. By the triangle and Rosenthal inequalities

$$\begin{aligned} \left\| \sum_{u=1}^{p_n} \sum_{k=-B_n}^{-l} \tilde{X}_u \tilde{X}_{u+k} \alpha_k \right\|_\tau &\le \sum_{h=1}^l \left\| \sum_{j=1}^{\lfloor (p_n-h)/l \rfloor} \sum_{k=-B_n}^{-l} \tilde{X}_{h+(j-1)l} \tilde{X}_{h+(j-1)l+k} \alpha_k \right\|_\tau \\ &\le O(l) \sqrt{p_n/l} \left\| \sum_{k=-B_n}^{-l} \tilde{X}_k \alpha_k \right\|_\tau \\ &\le O(\sqrt{p_n l}) \sum_{h=0}^{l-1} \left\| \sum_{j=0}^{\lfloor (B_n-l-h)/l \rfloor} \tilde{X}_{-B_n+h+jl} \alpha_{-B_n+h+jl} \right\|_\tau \\ &= O[(\sqrt{p_n l}) l \sqrt{B_n/l}] = O(\sqrt{p_n B_n} l). \end{aligned}$$

On the other hand, since $\tilde{X}_{h+3jl}\tilde{X}_{h+3jl+k}$, $0 \le j \le \lfloor (p_n - h)/(3l) \rfloor$, are i.i.d.,

$$
\left\| \sum_{u=1}^{p_n} \sum_{k=1-l}^{0} \tilde{X}_u \tilde{X}_{u+k}\alpha_k \right\|_\tau \le \sum_{k=1-l}^{0} \left\| \sum_{u=1}^{p_n} \tilde{X}_u \tilde{X}_{u+k}\alpha_k \right\|_\tau
$$

(A.13)
$$
\le \sum_{k=1-l}^{0} \sum_{h=1}^{3l} \left\| \sum_{j=0}^{\lfloor (p_n-h)/(3l) \rfloor} \tilde{X}_{h+3jl}\tilde{X}_{h+3jl+k}\alpha_k \right\|_\tau
$$

$$
= O(l^2 \sqrt{p_n/l}\,).
$$

Then we have $\|U_1\|_\tau = O(\sqrt{p_n B_n}\, l + l^2\sqrt{p_n/l}\,) = O(\sqrt{p_n B_n}\, l)$. Case (ii) $[B_n = O(\log n)]$. By the argument of (A.13), $\|U_1\|_\tau = O(B_n l\sqrt{p_n/l}\,) = O(\sqrt{p_n B_n}\, l)$. It is easily seen that the $O(\cdot)$-relation holds uniformly over $\lambda \in [0, \pi]$, that is,

(A.14)
$$
\sup_{\lambda \in [0,\pi]} \|U_1(\lambda)\|_\tau = O(l\sqrt{p_n B_n}\,).
$$

Then $\|U_1 - \mathbf{E}(U_1)\|_\tau = o[(nB_n)^{1/2}k_n^{-1/\tau}]$ and the Liapounov condition holds. By the central limit theorem and (A.11), we have (A.12). So (3.3) follows from (A.8). □

**A.3. Proof of Theorem 3.2.** We adopt the block method. Let $U_r(\lambda)$, $r = 1, \ldots, k_n$, be i.i.d. block sums with block length $p = p_n = \lfloor n^{1-4\eta/\delta}(\log n)^{-8/\delta-4} \rfloor$ and $V_r(\lambda)$, $r = 1, \ldots, k_n - 1$, be i.i.d. block sums with the same block length $q = q_n = p_n$. The last block $V_{k_n}(\lambda)$ is negligible. Note that $B_n = o(p_n)$ since $\eta < \delta/(4 + \delta)$. Let $l = l_n = \lfloor -8\log n/\log \rho(4) \rfloor$ as in the proof of Theorem 3.1. Define $U_r(\lambda)' := U_r(\lambda) \times \mathbf{1}(|U_r(\lambda)| \le d_n)$ for $r = 1, \ldots, k_n$, where $d_n = \lfloor \sqrt{nB_n}(\log n)^{-1/2} \rfloor$. The following lemma is needed.

LEMMA A.3. *Under the assumptions in Theorem 3.2, we have*

(A.15)
$$
\mathbf{E}\left( \max_{\lambda \in [0,\pi]} |V_{k_n}(\lambda)| \right) = O(\sqrt{p_n l}\, B_n),
$$

(A.16)
$$
\mathbf{E}\left( \max_{\lambda \in [0,\pi]} |h_n(\lambda)| \right) = o(1),
$$

(A.17)
$$
\max_{\lambda \in [0,\pi]} \mathrm{var}(U_1(\lambda)) = O(p_n B_n),
$$

(A.18)
$$
\mathrm{var}(U_1(\lambda)') = \mathrm{var}(U_1(\lambda))[1 + o(1)],
$$

*where the relation $o(1)$ in (A.18) holds uniformly over $[0, \pi]$.*

PROOF. Let $z = k_n(p+q) + 1 - q$ and $\tau = 2 + \delta/2$. Then

$$\mathbf{E}\left(\max_{\lambda \in [0,\pi]} |V_{k_n}(\lambda)|\right) \leq C \sum_{j=-B_n}^{B_n} \mathbf{E}\left|\sum_{u=z}^{n} \tilde{X}_u \tilde{X}_{u+j}\right|.$$

For $|j| \leq l$, $\|\sum_{u=z}^{n} \tilde{X}_u \tilde{X}_{u+j}\| = O(\sqrt{p_n l})$ since $\tilde{X}_u \tilde{X}_{u+j}$ is $2l$-dependent. When $|j| > l$, $\|\sum_{u=z}^{n} \tilde{X}_u \tilde{X}_{u+j}\|^2 = \sum_{u,u'=z}^{n} \mathbf{E}(\tilde{X}_u \tilde{X}_{u+j} \tilde{X}_{u'} \tilde{X}_{u'+j}) = O(p_n l)$ since the sum vanishes if $|u - u'| > l$. So $\mathbf{E} \max_{\lambda \in [0,\pi]} |V_{k_n}(\lambda)| = O(\sqrt{p_n l} B_n)$. Let $\tilde{h}_n(\lambda)$ be the corresponding sum of $h_n(\lambda)$ with $X_u X_{u+k}$ replaced by $\tilde{X}_u \tilde{X}_{u+k}$. As at (A.9), we have $\mathbf{E} \max_{\lambda \in [0,\pi]} |h_n(\lambda) - \tilde{h}_n(\lambda)| = o(1)$. To show (A.16), it suffices to show $\mathbf{E} \max_{\lambda \in [0,\pi]} |\tilde{h}_n(\lambda)| = o(1)$ which follows from a similar argument as in the proof of (A.15). Regarding (A.17), we have

$$\begin{aligned}
\text{var}(U_1(\lambda)) &= \left\|\sum_{u=1}^{p} \sum_{k=-B_n}^{B_n} \{X_u X_{u+k} - r(k)\}\alpha_k\right\|^2 \\
&= \sum_{u,u'=1}^{p} \sum_{k,k'=-B_n}^{B_n} \{r(u-u')r(u-u'+k-k') \\
&\qquad + r(u'-u+k')r(u'-u-k) \\
&\qquad + \text{cum}(X_0, X_k, X_{u'-u}, X_{u'-u+k'})\}\alpha_k \alpha_{k'} \\
&=: I_1 + I_2 + I_3.
\end{aligned}$$

Then $I_1$ is bounded by $C \sum_{h=1-p}^{p-1} (p - |h|)|r(h)| \sum_{g=-2B_n}^{2B_n} (2B_n + 1 - |g|)|r(h+g)|$, which is less than $Cp(2B_n+1)(\sum_{k=-\infty}^{\infty} |r(k)|)^2$. Similarly, smaller bounds can be obtained for $I_2$ and $I_3$ due to the summability of the second and fourth cumulants. Thus $\max_{\lambda \in [0,\pi]} \text{var}(U_1(\lambda)) = O(p_n B_n)$. For (A.18), let $v = \text{var}\{U_1(\lambda) - U_1(\lambda)'\}$ and $c = \mathbf{E}(U_1(\lambda)')\mathbf{E}\{U_1(\lambda) - U_1(\lambda)'\}$. Then $\text{var}(U_1(\lambda)') = \text{var}(U_1(\lambda)) - v + 2c$. By Markov's inequality and (A.14), $v \leq \|U_1(\lambda)\|_{\tau}^{\tau}/d_n^{\tau-2} = o(p_n B_n)$ and similarly $c \leq \|U_1(\lambda)\|_{\tau}^{\tau+1}/d_n^{\tau-1} = o(p_n B_n)$. By Lemma A.2 and since $f$ is everywhere positive, (A.18) follows. $\square$

PROOF OF THEOREM 3.2. Let $H_n(\lambda) = \sum_{r=1}^{k_n} [U_r(\lambda) - \mathbf{E}\{U_r(\lambda)\}]$, $H_n(\lambda)' = \sum_{r=1}^{k_n} [U_r(\lambda)' - \mathbf{E}\{U_r(\lambda)'\}]$. Let $\lambda_j = \pi j/t_n, j = 0, \ldots, t_n$, $t_n = \lfloor B_n \log(B_n) \rfloor$. Let $c_n = 1/(1 - 3\pi/\log B_n) \to 1$. By Corollary 2.1 in [72], $\max_{\lambda \in [0,\pi]} |H_n(\lambda)| \leq c_n \max_{j \leq t_n} |H_n(\lambda_j)|$. By (A.17) and (A.18), there exists a constant $C_1 > 1$ such that $\max_{\lambda \in [0,\pi]} \text{var}(U_1(\lambda)') \leq C_1 p_n B_n$. Let $\alpha_n = (C_1 n B_n \log n)^{1/2}$. By Bernstein's inequality, we have

$$\begin{aligned}
\mathbf{P}\left(\max_{0 \leq j \leq t_n} |H_n(\lambda_j)'| \geq 4\alpha_n\right) &\leq \sum_{j=0}^{t_n} \mathbf{P}(|H_n(\lambda_j)'| \geq 4\alpha_n) \\
&= O(t_n)\exp\left(\frac{-16\alpha_n^2}{2k_n C_1 p_n B_n + 16 d_n \alpha_n}\right) = o(1).
\end{aligned}$$

Let $U_r(\lambda)'' = U_r(\lambda) - U_r(\lambda)'$ and $H_n(\lambda)'' = H_n(\lambda) - H_n(\lambda)'$. Then by Markov's inequality and (A.14),

$$
\mathbf{P}\Big(\max_{0 \le j \le t_n} |H_n(\lambda_j)''| \ge 4\alpha_n\Big) \le \sum_{j=0}^{t_n} \mathbf{P}(|H_n(\lambda_j)''| \ge 4\alpha_n)
$$

$$
\le \sum_{j=0}^{t_n} \frac{\operatorname{var}(U_1(\lambda_j)'')k_n}{16\alpha_n^2} = \frac{O(t_n k_n (\sqrt{p_n B_n} l)^\tau)}{\alpha_n^2 d_n^{\tau-2}}
$$

$$
= O\Big(\frac{(B_n \log n)(n/p_n)(\sqrt{p_n B_n} \log n)^\tau}{(nB_n \log n)(nB_n)^{\tau/2-1}(\log n)^{-\tau/2+1}}\Big)
$$

$$
= O((\log n)^{-\delta/4}) = o(1).
$$

So $\max_{\lambda \in [0,\pi]} |H_n(\lambda)| = O_{\mathbf{P}}(\alpha_n)$. Clearly the same bound also holds for the sum $\sum_{r=1}^{k_n-1} [V_r(\lambda) - \mathbf{E}\{V_r(\lambda)\}]$. By (A.9), $\mathbf{E}\max_{\lambda \in [0,\pi]} |\tilde{g}_n(\lambda) - g_n(\lambda)| = o(1)$. By (A.15), (A.16) and (A.8), we have (3.4). $\square$

**A.4. Proof of Propositions 4.1 and 4.2.** Let $\delta_{j,k} = \mathbf{1}_{j=k}$.

LEMMA A.4.  *Let $m = \lfloor (n-1)/2 \rfloor$. (i) Assume* (A4'), (A5') *and $X_t \in \mathcal{L}^8$. Then $\max_{j,k \le m} |\operatorname{cov}(I_j^2, I_k^2) - 4f_j^4 \delta_{j,k}| = O(1/n)$. (ii) Assume* (A4'), (A5) *and $X_t \in \mathcal{L}^4$. Then $\max_{j,k \le m} |\operatorname{cov}(I_j, I_k) - f_j^2 \delta_{j,k}| = O(1/n)$.*

PROOF.  We only show (i) since (ii) can be handled similarly. Note that

$$
\operatorname{cov}(I_j^2, I_k^2)
$$

$$
(A.19) \quad = \frac{1}{16\pi^4 n^4} \sum_{t_l, s_l \in \{1,\dots,n\}, l=1,\dots,4} e^{i(t_1-t_2+t_3-t_4)\lambda_j - i(s_1-s_2+s_3-s_4)\lambda_k}
$$

$$
\times \operatorname{cov}(X_{t_1} X_{t_2} X_{t_3} X_{t_4}, X_{s_1} X_{s_2} X_{s_3} X_{s_4}).
$$

By Theorem II.2 in [57], we have

$$
\operatorname{cov}(X_{t_1} X_{t_2} X_{t_3} X_{t_4}, X_{s_1} X_{s_2} X_{s_3} X_{s_4}) = \sum_v \operatorname{cum}(X_{i_j}; i_j \in v_1) \cdots \operatorname{cum}(X_{i_j}; i_j \in v_p),
$$

where $\sum_v$ is over all indecomposable partitions $v = v_1 \cup \cdots \cup v_p$ of the two-way table

$$
\begin{array}{llll}
X_{t_1}(+) & X_{t_2}(-) & X_{t_3}(+) & X_{t_4}(-) \\
X_{s_1}(-) & X_{s_2}(+) & X_{s_3}(-) & X_{s_4}(+).
\end{array}
$$

The signs in the above table are from the exponential terms in the sum (A.19). Since $\mathbf{E}(X_t) = 0$, only partitions $v$ with $\#v_j > 1$ for all $j$ contribute.

One of the many indecomposable partitions consisting only of pairs with $+$ in $t$ matched to $-$ in $s$ [say, $\{(t_1, s_1), (t_2, s_2), (t_3, s_3), (t_4, s_4)\}$] leads to the sum $[A(\lambda_j, \lambda_k)]^4$, where

$$A(\lambda_j, \lambda_k) = \frac{1}{2\pi n} \sum_{t_1, s_1 = 1}^{n} r(t_1 - s_1)e^{it_1\lambda_j - is_1\lambda_k} = f(\lambda_j)\mathbf{1}_{j=k} + O(1/n).$$

The other indecomposable partitions consisting entirely of pairs (with $+$ in $t$ matched to $-$ in $s$) are $\{(t_1, s_3), (t_2, s_2), (t_3, s_1), (t_4, s_4)\}$, $\{(t_1, s_1), (t_2, s_4), (t_3, s_3), (t_4, s_2)\}$ and $\{(t_1, s_3), (t_2, s_4), (t_3, s_1), (t_4, s_2)\}$. It is easily seen after some calculations that partitions containing entirely pairs but with at least one $+$ in $t$ matched to one $+$ in $s$ result in a term of order $O(1/n)$ for any $j, k$. All other partitions that are not all pairs will give a quantity of order $O(1/n)$ due to the summability of cumulants up to the eighth order. Finally, it is not hard to see that $O(1/n)$ does not depend on $(j, k)$. Thus the conclusion is proved. $\square$

LEMMA A.5. *Assume $X_t \in \mathcal{L}^8$, (A2), (A3$'$), (A4$'$) and (A5$'$). Then* $\mathrm{var}^*(\varepsilon_1^*) \to 1$ *in probability and* $\mathbf{E}^*(|\varepsilon_1^*|^4) = O_{\mathbf{P}}(1)$.

PROOF. By (A3$'$), $\tilde{f}$ is a uniformly consistent estimate of $f$. It remains to show

$$(A.20) \quad \frac{1}{N}\sum_{j=1}^{N}\frac{I_j}{f_j} \to 1, \qquad \frac{1}{N}\sum_{j=1}^{N}\frac{I_j^2}{f_j^2} \to 2 \qquad \text{in probability and}$$

$$\frac{1}{N}\sum_{j=1}^{N}\frac{I_j^4}{f_j^4} = O_{\mathbf{P}}(1).$$

By Proposition 10.3.1 in [8] and Lemma A.4, we have $\mathbf{E}(I_j) = f_j + o(1)$ and $\mathbf{E}(I_j^2) = 2f_j^2 + o(1)$ uniformly in $j$. Thus the first two assertions follow from Lemma A.4 since their variances go to 0 as $n \to \infty$. By Lemma A.4, $\mathbf{E}(I_j^4) = \mathrm{cov}(I_j^2, I_j^2) + (\mathbf{E}I_j^2)^2 = 8f_j^4 + o(1)$ uniformly in $j$, and the last assertion holds. $\square$

REMARK A.1. For linear processes, Franke and Härdle [24] remarked that their consistency result strongly depends on the asymptotic normality of $f_n$ and the weak convergence of $F_{\tilde{I}, m}(x)$ (see Corollary 2.2). The latter condition holds under $\varepsilon_1 \in \mathcal{L}^5$ and (4.2) by Chen and Hannan [12]. Franke and Härdle [24] further conjectured that their results hold assuming only $\varepsilon_1 \in \mathcal{L}^4$, under which the weak convergence of $F_{\tilde{I}, m}(x)$ might be true. However, it seems from our argument (see the proof of Proposition 4.1) that it is not the weak convergence of $F_{\tilde{I}, m}(x)$ but the first two conditions in (A.20) that play key roles; compare Proposition A1 in [24]. The proof of the second assertion in (A.20) (see Lemmas A.4 and A.5) in a general setting needs the stronger eighth moment assumption.

Let $a \vee b = \max(a, b)$ and $a \wedge b = \min(a, b)$; let $\tilde{r}_2(k) = \int_0^{2\pi} \tilde{f}^2(\lambda) e^{ik\lambda} \, d\lambda$, $r_2(k) = \int_0^{2\pi} f^2(\lambda) e^{ik\lambda} \, d\lambda$, $\tilde{r}(k) = \int_0^{2\pi} \tilde{f}(\lambda) e^{ik\lambda} \, d\lambda$ and $F_n^+ = \{1, \ldots, \lfloor n/2 \rfloor\}$. By (A3), $\max_{k \in \mathbf{Z}} |\tilde{r}_2(k) - r_2(k)| \leq 2\pi \max_\lambda |\tilde{f}^2(\lambda) - f^2(\lambda)| = o_{\mathbf{P}}(b_n)$.

PROOF OF PROPOSITION 4.1. By Lemma 8.3 of [4], the convergence under the $d_2$ metric is equivalent to weak convergence and convergence of the first two moments. By (A6), it suffices to show that

$$(A.21) \quad nb_n \operatorname{var}^*(f_n^*(\lambda)) \to \sigma^2(\lambda), \qquad \mathcal{L}(V_n^*(\lambda)|\{X_j\}_{j=1}^n) \Rightarrow N(0, \sigma^2(\lambda))$$

$$\text{in probability.}$$

Let $\Delta_j = \sum_{k=-B_n}^{B_n} a(kb_n) e^{-ik\lambda}(e^{ik\omega_j} + e^{-ik\omega_j})$. Since the resampled residuals $\{\varepsilon_j^*\}$ are i.i.d. given $X_1, \ldots, X_n$, we have $\operatorname{var}^*(I_j^*) = \tilde{f}_j^2 \operatorname{var}^*(\varepsilon_1^*)$, and, since $I_0^* = 0$, $nb_n \operatorname{var}^*(f_n^*(\lambda)) = \operatorname{var}^*(\varepsilon_1^*) R_n(\lambda) + o_{\mathbf{P}}(1)$, where

$$\begin{aligned}
R_n(\lambda) &= \frac{nb_n}{n^2} \sum_{j \in F_n^+} \tilde{f}_j^2 \Delta_j^2 \\
&= \frac{1}{nB_n} \sum_{k,k'=-B_n}^{B_n} a(kb_n) a(k'b_n) e^{-i\lambda(k-k')} \\
&\qquad \times \sum_{j \in F_n} \tilde{f}_j^2 \{e^{i\omega_j(k-k')} + e^{i\omega_j(k+k')}\} + o_{\mathbf{P}}(1) \\
&= \frac{1}{2\pi B_n} \sum_{k,k'=-B_n}^{B_n} a(kb_n) a(k'b_n) e^{-i\lambda(k-k')} \\
&\qquad \times \{\tilde{r}_2(k-k') + \tilde{r}_2(k+k')\} + o_{\mathbf{P}}(1) \\
&= \frac{1}{2\pi B_n} \sum_{k,k'=-B_n}^{B_n} a(kb_n) a(k'b_n) e^{-i\lambda(k-k')} \\
&\qquad \times \{r_2(k-k') + r_2(k+k')\} + o_{\mathbf{P}}(1) \\
&= R_n^{(1)}(\lambda) + R_n^{(2)}(\lambda) + o_{\mathbf{P}}(1) \qquad \text{(say)}.
\end{aligned}$$

Let $\beta_n(k) = \int_0^{2\pi} R_n^{(1)}(\lambda) e^{ik\lambda} \, d\lambda$ and $\beta(k) = \int_0^{2\pi} \int_{-1}^1 a^2(u) f^2(\lambda) e^{ik\lambda} \, du \, d\lambda$. Then

$$\beta_n(k) = \frac{r_2(k)}{B_n} \sum_{j=-B_n+0 \vee k}^{B_n+0 \wedge k} a(jb_n) a((j-k)b_n) \to r_2(k) \int_{-1}^1 a^2(u) \, du.$$

Since $|\beta_n(k)| \leq C|r_2(k)|$ and $\sum_{k \in \mathbf{Z}} |r_2(k)| < \infty$, by the Lebesgue dominated convergence theorem, $R_n^{(1)}(\lambda) \to f^2(\lambda) \int_{-1}^{1} a^2(u)\,du$. For $R_n^{(2)}(\lambda)$, $\lambda \neq 0, \pm\pi$, we have

$$R_n^{(2)}(\lambda) = \frac{1}{2\pi B_n} \sum_{h=-2B_n}^{2B_n} r_2(h)e^{ih\lambda} \sum_{k=-B_n+0\vee h}^{B_n+0\wedge h} a(kb_n)a((k-h)b_n)e^{-2ik\lambda}$$

$$= \frac{1}{2\pi B_n} \sum_{h=-2B_n}^{2B_n} r_2(h)e^{ih\lambda}O(1) \to 0.$$

It is easily seen that $R_n^{(1)}(\lambda) = R_n^{(2)}(\lambda)$ when $\lambda = 0, \pm\pi$. Hence by Lemma A.5, $nb_n \operatorname{var}^*(f_n^*(\lambda)) \to \sigma^2(\lambda)$ in probability.

Finally, since $\{\varepsilon_j^*\}$ are i.i.d. conditional on $\{X_1, \ldots, X_n\}$, by the Berry–Esseen theorem and Lemma A.5, we have

$$\sup_x \left| \mathbf{P}^*(V_n^*(\lambda) \leq x) - \Phi\left(\frac{x}{nb_n \operatorname{var}^*(f_n^*(\lambda))}\right) \right|$$

$$\leq \frac{C \sum_{j \in F_n^+} \tilde{f}_j^4 \mathbf{E}^* |\varepsilon_1^*|^4 \Delta_j^4}{[\sum_{j \in F_n^+} \tilde{f}_j^2 \operatorname{var}^*(\varepsilon_1^*)\Delta_j^2]^2}$$

$$= O_{\mathbf{P}}\left(\frac{nB_n^4}{n^2 B_n^2}\right),$$

which implies $\mathcal{L}(V_n^*(\lambda)|X_1, \ldots, X_n) \Rightarrow N(0, \sigma^2(\lambda))$ in probability since $B_n^2 = o(n)$ and $\sup_x |\Phi(x/\sigma_1) - \Phi(x/\sigma_2)| \leq C|\sigma_1/\sigma_2 - 1|$ for some constant $C$. □

PROOF OF PROPOSITION 4.2. Since $B_n^3 = o(n)$ and $\tilde{r}(k) = a(k\tilde{b}_n)\hat{r}(k)$, $|k| \leq \tilde{B}_n$ and 0 otherwise, we have $B_n^2[\mathbf{E}^* f_n^*(\lambda) - \tilde{f}(\lambda)] = J_n(\lambda) + o_{\mathbf{P}}(1)$, where

$$J_n(\lambda) = \frac{B_n^2}{2\pi} \sum_{k=-\tilde{B}_n}^{\tilde{B}_n} a(k\tilde{b}_n)\hat{r}(k)e^{-ik\lambda}(a(kb_n) - 1).$$

It remains to show $\mathbf{E}(J_n(\lambda)) \to c_2 f''(\lambda)$ and $\operatorname{var}(J_n(\lambda)) \to 0$. By (A1), (A4) and (A5),

$$\mathbf{E}(J_n(\lambda)) = \frac{B_n^2}{2\pi} \sum_{k=-\tilde{B}_n}^{\tilde{B}_n} a(k\tilde{b}_n)e^{-ik\lambda}(1 - |k|/n)r(k)(a(kb_n) - 1)$$

$$= -\frac{B_n^2}{2\pi} \sum_{k=-\tilde{B}_n}^{\tilde{B}_n} a(k\tilde{b}_n)e^{-ik\lambda}r(k)k^2 b_n^2 c_2(1 + o(1)) \to c_2 f''(\lambda)$$

and

$$\mathrm{var}(J_n(\lambda)) = \frac{B_n^4}{4\pi^2} \sum_{k,k'=-\tilde{B}_n}^{\tilde{B}_n} a(k\tilde{b}_n)a(k'\tilde{b}_n)(a(kb_n)-1)(a(k'b_n)-1)$$

$$\times e^{-i(k-k')\lambda}\,\mathrm{cov}(\hat{r}(k),\hat{r}(k'))$$

$$(A.22)\qquad = \frac{(1+o(1))c_2^2}{4\pi^2 n^2} \sum_{k,k'=-\tilde{B}_n}^{\tilde{B}_n} a(k\tilde{b}_n)a(k'\tilde{b}_n)k^2 k'^2 e^{-i(k-k')\lambda}$$

$$\times \sum_{t=1}^{n-|k|}\sum_{t'=1}^{n-|k'|}\mathrm{cov}(X_t X_{t+|k|}, X_{t'} X_{t'+|k'|})$$

$$= O(\tilde{B}_n^4/n^2) \sum_{k,k'=0}^{\tilde{B}_n}\sum_{t=1}^{n-k}\sum_{t'=1}^{n-k'}|\mathrm{cov}(X_t X_{t+k}, X_{t'} X_{t'+k'})|.$$

Note that $\mathrm{cov}(X_t X_{t+k}, X_{t'} X_{t'+k'}) = r(t-t')r(t+k-t'-k') + r(t-t'-k')r(t'-t-k) + \mathrm{cum}(X_t, X_{t+k}, X_{t'}, X_{t'+k'})$. The contribution of the first term $r(t-t')r(t+k-t'-k')$ to (A.22) is $O(\tilde{B}_n^5/n)\sum_{h=-\tilde{B}_n}^{\tilde{B}_n}\sum_{s=-2n}^{2n}|r(h)r(h+s)| = O(\tilde{B}_n^5/n) = o(1)$ since $\sum_{k\in\mathbf{Z}}|r(k)| < \infty$. Similarly, the contribution of the second term to (A.22) approaches zero as $n \to \infty$. The third term is $O(\tilde{B}_n^4/n) = o(1)$ due to the summability of the fourth cumulants. $\square$

**Acknowledgments.** The authors would like to thank Michael Stein for helpful comments on an earlier version. Thanks are also due to the referees for their thoughtful and constructive suggestions.

Department of Statistics
University of Illinois at Urbana-Champaign
725 South Wright Street
Champaign, Illinois 61820
USA
E-mail: xshao@uiuc.edu

Department of Statistics
University of Chicago
5734 S. University Avenue
Chicago, Illinois 60637
USA
E-mail: wbwu@galton.uchicago.edu